\documentclass[11pt,a4paper]{article}
\usepackage[intlimits]{amsmath}
\usepackage{graphicx}
\usepackage{amsthm}
\usepackage{amssymb}
\newcommand{\CC}{\mathcal{C}}
\newcommand{\del}{\delta}
\newcommand{\eps}{\varepsilon}
\newcommand{\lam}{\lambda}
\newcommand{\sig}{\sigma}

\newcommand{\Del}{\Delta}
\newcommand{\Ome}{\Omega}
\newcommand{\dO}{{\partial\Ome}}
\newcommand{\Sbold}{\mathbf{S}}
\newcommand{\nbold}{\mathbf{n}}
\newcommand{\xbold}{\mathbf{x}}
\newcommand{\vbold}{\mathbf{v}}
\newcommand{\Rbb}[1][]{\mathbb{R}^{#1}}
\newcommand{\dsp}{\displaystyle}
\begin{document}
\numberwithin{equation}{section}
\theoremstyle{plain}
\newtheorem{thm}{Theorem}[section]
\newtheorem{lem}[thm]{Lemma}
\newtheorem{prop}[thm]{Proposition}
\newtheorem{cor}[thm]{Corollary}
\newtheorem{conj}[thm]{Conjecture}
\theoremstyle{definition}
\newtheorem{defn}[thm]{Definition}
\newtheorem{exmp}[thm]{Example}
\theoremstyle{remark}
\newtheorem{rem}[thm]{Remark}
\newtheorem{note}[thm]{Note}
\title{Range of the first three eigenvalues of the planar Dirichlet Laplacian}
\author{Michael Levitin\thanks{The research of M.L. was
partially supported by EPSRC grant GR/M20990 and by EPSRC Spectral Theory Network}\\
Rustem Yagudin\thanks{The research of R.Y. was
partially supported by ORS grant}\\
\ \\
\normalsize\small Department of Mathematics, Heriot-Watt University\\
\normalsize\small Riccarton, Edinburgh EH14 4AS, U.~K.\\
\normalsize\small email {\sffamily M.Levitin@ma.hw.ac.uk}\\
\normalsize\small www {\sffamily http://www.ma.hw.ac.uk/$\sim$levitin}%
}
\date{}
\maketitle
\begin{abstract} We conduct extensive numerical experiments aimed at 
finding the admissible range of the ratios of the first three eigenvalues of
a planar Dirichlet Laplacian. The results improve the previously known
theoretical estimates of M~Ashbaugh and R~Benguria. We also prove some properties of a maximizer of the ratio $\lambda_3/\lambda_1$. 
\end{abstract}
\section{Introduction}

Let $\Ome$ be a bounded domain in $\Rbb[n]$, $n\ge 2$. We consider the
eigenvalue problem for the Dirichlet Laplacian,
\begin{equation}\label{eq:Lapl}
-\Del u=\lam u\qquad\text{in\quad $\Ome$}\,,
\end{equation}
\begin{equation}\label{eq:Dir}
u|_\dO=0\,.
\end{equation}

Let us denote the eigenvalues by $\lam_1(\Ome)$, $\lam_2(\Ome)$,  \dots, (we will
sometimes omit explicit dependence on $\Ome$ when speaking about a generic domain),
where $0<\lam_1<\lam_2\le\lam_3\le\dots$. The corresponding
orthonormal basis of real eigenfunctions will be denoted $\dsp\{u_j\}_{j=1}^\infty$.

For the last fifty years, the problem of obtaining {\it a priori\/} estimates
of the eigenvalues and their ratios has attracted a substantial attention.
The existing results can be roughly divided into two groups --- {\it universal\/} estimates, valid, as the name suggests, for all eigenvalues and
all the domains in $\Rbb[n]$, which do not take into account any geometric information, and {\it isoperimetric\/} estimates for low eigenvalues. We
briefly survey some known results below; the reader is referred to the very detailed survey
paper \cite{Ash} and references therein for a full discussion.

\subsection*{Universal estimates} Probably the first, and best known, estimate
of this type is the Payne-P\'olya-Weinberger inequality \cite{PPW},
\begin{equation*}\label{eq:PPW}
\lam_{m+1}\le \lam_m+\frac{4}{mn}\sum_{j=1}^m\lam_j\,.\tag{PPW}
\end{equation*}
This was subsequently improved by Hile and Protter \cite{HiPr}, and, in 1990's,
by Hong Cang Yang \cite{HCY}, whose implicit estimate
\begin{equation*}\label{eq:HCY}
\sum_{j=1}^m \left(\lam_{m+1}-\lam_j\right)\left(\lam_{m+1}-\left(1+\frac{4}{n}\right)
\lam_j\right)\le 0\,,\tag{HCY}
\end{equation*}
remains the best universal estimate so far for the eigenvalues of the Dirichlet Laplacian.

The general method of obtaining \eqref{eq:PPW} and \eqref{eq:HCY}, as well as of similar estimates for
a variety of other operators, has been  the use of variational principles with some ingenious choices
of trial functions, see \cite{Ash}. Recently, an alternative abstract scheme, based on the so-called commutator
trace identities, which easily implies, in particular, \eqref{eq:PPW} and \eqref{eq:HCY}, has been developed in \cite{LevPar}, see also \cite{HarStu}.

By their very nature, the universal estimates are generically non-sharp.

\subsection*{Isoperimetric estimates} Both \eqref{eq:PPW} and \eqref{eq:HCY} give, for $m=1$, the
estimate
$$
\frac{\lam_2}{\lam_1}\le 1+\frac{4}{n}\,.
$$
This upper bound cannot, in fact, be attained. Already, Payne, P\'olya, and Weinberger conjectured
that the actual optimal upper bound on the ratio of the first two eigenvalues of the Dirichlet
Laplacian is
\begin{equation*}\label{eq:Ks}
\frac{\lambda_2}{\lambda_1}(\Ome)\le
\left.\frac{\lambda_2}{\lambda_1}\right|_{n\text{-dimensional ball}}=
\frac{j_{n/2,1}^2}{j_{n/2-1,1}^2}=:K_n\qquad\text{for }\Ome\subset\Rbb[n]
\tag{AB${}_0$}
\end{equation*}
(here $j_{p,q}$ stands for the $q$-th zero of the Bessel function $J_p(\rho)$, so, in a planar case
$n=2$,   $K_2\approx 2.5387$ compared with \eqref{eq:PPW} bound
$\dsp \left.\frac{\lambda_2}{\lambda_1}(\Ome)\right|_{\Ome\subset\Rbb[2]}\le 3$).

Conjecture \eqref{eq:Ks} was eventually proved, only in the early 1990s, by Ashbaugh and Benguria \cite{AshBen1, AshBen3},
using, in particular,  symmetrization techniques going back to
the Faber-Krahn inequality,
$$
\lambda_1(\Omega)\ge \lambda_1(\Omega^\star)\,,
$$
where $\Omega^\star$ is an $n$-dimensional ball of the same volume as $\Ome$.

We would like to mention, in this context, extensive computational experiments designed to verify
\eqref{eq:Ks} by Haeberly \cite{Hae, HaeOve}.

\subsection*{Statement of the problem} As mentioned, \eqref{eq:Ks} gives the full description
of the range of the possible values of ratio of the first two eigenvalues of the Dirichlet Laplacian,
$\dsp \frac{\lambda_2}{\lambda_1}$ for domains in the
Euclidean space (the obvious lower bound is $\dsp \frac{\lambda_2}{\lambda_1}\ge 1$).
In fact, similar results were also obtained for domains in $\mathbb{S}^n$ and
$\mathbb{H}^n$. A natural extension would be to find optimal upper bounds on the range of
the ratios of the first {\it three\/} eigenvalues of the Dirichlet Laplacian,
$\dsp \left(\frac{\lambda_2}{\lambda_1}, \frac{\lambda_3}{\lambda_1}\right)$, in particular
for planar domains. In other words, we would like to find, for
$\dsp x:=\frac{\lambda_2}{\lambda_1}$ and
$\dsp y:=\frac{\lambda_3}{\lambda_1}$, the function
\begin{equation}\label{eq:yast}
y^\ast(x):=\max_{\Omega\subset\Rbb[2]:\frac{\lambda_2}{\lambda_1}(\Ome)=x}
\frac{\lambda_3}{\lambda_1}(\Ome)
\end{equation}
and the number
\begin{equation}\label{eq:Yast}
Y^\ast:=\max_{x\in[1,K_2]} y^\ast(x)=\max_{\Omega\subset\Rbb[2]} \frac{\lambda_3}{\lambda_1}(\Ome)
\end{equation}
or their best possible estimates. We will use notation \eqref{eq:yast} and \eqref{eq:Yast} when
looking for maxima in particular classes of domains as well.

Despite an apparent simplicity of this problem, and a wide attention
it has attracted, it turned out to be rather difficult. In \cite{AshBen5, AshBen6}, Ashbaugh and Benguria proved a complicated 
upper bound for $y^\ast(x)$ and also demonstrated that
\begin{equation}\label{eq:lam3best}
3.1818\approx \frac{35}{11}\le Y^\ast \lessapprox 3.83103\,.
\end{equation}
Their estimates improve upon previous results due to  Payne, P\'olya, and Weinberger, Brands, de Vries,
Hile and Protter, Marcellini, Chiti, Hong Cang Yang, and themselves; see \cite{AshBen5, AshBen6} and their earlier papers  \cite{AshBen2, AshBen4} for extensive
bibliography and details of proofs. We present their estimates and other known facts in the next Section; just note at the moment that the lower bound in \eqref{eq:lam3best} is attained when $\Ome$ is the rectangle
$R_a:=[0,1]\times[0,a]$ with $\dsp a=\sqrt\frac{8}{3}$.

In the current paper, we describe extensive numerical experiments aimed at improving \eqref{eq:lam3best}. We also show, using perturbation techniques, that  the rectangle
$\dsp R_{\sqrt\frac{8}{3}}$ does not maximize the ratio $\dsp \frac{\lambda_3}{\lambda_1}$ and indicate a class of domains among which a possible maximizer could be
found.

\subsection*{Acknowledgements} The original suggestion to conduct numerical experiments on the range
of  $\dsp \left(\frac{\lambda_2}{\lambda_1}, \frac{\lambda_3}{\lambda_1}\right)$ for planar domains came from 
Brian Davies; we would like to thank him, as well as Mark Ashbaugh and Leonid
Parnovski, for valuable discussions and advice. It is the authors, nevertheless, who take full
responsibility for the realization of this idea, and any criticism for possible shortcomings of this
realization should be addressed to them.

\section{Known results for the range of $\dsp \left(\frac{\lambda_2}{\lambda_1}, \frac{\lambda_3}{\lambda_1}\right)$ for planar domains}

\subsection*{Explicit solutions} The spectral problem \eqref{eq:Lapl}, \eqref{eq:Dir} admits a full
solution by separation of variables when $\Ome$ is, for example, a disjoint union of a number of rectangles or circles.  For a reference, we collect below the results on the range of 
$\dsp \left(\frac{\lambda_2}{\lambda_1}, \frac{\lambda_3}{\lambda_1}\right)$ in these cases.

\begin{description}
\item[Rectangles] Let  $R_a:=[0,1]\times[0,a]$ be a rectangle with the side ratio $a$; without
loss of generality $a\ge 1$. Then
$$
\frac{\lambda_2}{\lambda_1}(R_a)=\frac{a^2+4}{a^2+1}
$$
and
$$
\dfrac{\lambda_3}{\lambda_1}(R_a)=\begin{cases}
\dfrac{a^2+9}{a^2+1}\qquad\text{for}\quad a \geq \sqrt{\dfrac{8}{3}}\,,\\
\dfrac{4a^2+1}{a^2+1}\qquad\text{for}\quad 1\leq a \leq \sqrt{\dfrac{8}{3}}\,.
\end{cases}
$$
Thus, for rectangles, in notation \eqref{eq:yast} and \eqref{eq:Yast},
\begin{equation}\label{eq:yastrect}
\left.y(x)=y^\ast(x)\right|_{\text{rectangles}}=\begin{cases}
\dfrac{8}{3}x-\dfrac{5}{3}\qquad\text{for}\quad 1\leq x \leq \dfrac{20}{11}\,,\\
5-x\qquad\text{for}\quad \dfrac{20}{11}\leq x \leq \dfrac{5}{2}\,,
\end{cases}
\end{equation}
and the maximum value of $\dsp \frac{\lambda_3}{\lambda_1}$ is
$$
\left.Y^\ast\right|_{\text{rectangles}}=\dfrac{35}{11}\,,
$$
attained when $\dsp a=\sqrt\frac{8}{3}$. Note that for this particular rectangle $\lam_3$ is
a degenerate eigenvalue: $\dsp \lam_3(R_{\sqrt\frac{8}{3}})=\lam_4(R_{\sqrt\frac{8}{3}})$,
and it is the only $a$ for which $\lam_3(R_a)$ is not simple.

In the $(x,y)$-plane, \eqref{eq:yastrect} corresponds to the two straight lines intersecting at the
point $\dsp \left(\frac{20}{11},\frac{35}{11}\right)$.
\item[Circles] For a single circle, $x=y=K_2$. As easily checked, for a union of more than one
disjoint circles of arbitrary radii,
\begin{equation}\label{eq:yastcirc}
\left.y^\ast(x)\right|_{\text{circles}}\equiv K_2\,,\qquad\text{for}\quad 1\le x\le K_2\,.
\end{equation}
Its graph in the $(x,y)$-plane is a straight line parallel to the $x$-axis.
\item[Disjoint unions] The following easily checked fact shows that one cannot obtain higher values of
$y^\ast(x)$ by considering disjoint unions of sets from two different classes. Namely, let $\CC_j$ be two arbitrary
classes of domains, with corresponding functions $\left.y^\ast(x)\right|_{\CC_j}$ (not necessarily defined for all $x\in[1,K_2]$).
Then, for any domain $\Omega=\Omega_1\sqcup\Omega_2$ with $\Omega_j\in\CC_j$, we have, for $\dsp x=\frac{\lambda_2}{\lambda_1}(\Omega)$ and
$\dsp y=\frac{\lambda_3}{\lambda_1}(\Omega)$, the inequality $y\le\max\left(\left.y^\ast(x)\right|_{\CC_1}, \left.y^\ast(x)\right|_{\CC_2}, K_2\right)$.
\item[Other domains] There are other domains, like sectors of the annuli, ellipses, etc., for which the
problem of funding the eigenvalues is reduced by separation of variables to the problem of solving some transcendental
equations. However, the latter one is often not easier than the numerical solution of the original problem,
so we do not treat these cases here.
\end{description}

The graphs of $\dsp \left.y(x)=y^\ast(x)\right|_{\text{rectangles}}$ and
$\dsp \left.y^\ast(x)\right|_{\text{circles}}$ are shown in Figure~\ref{fig:figure_1}.

\subsection*{Ashbaugh-Benguria estimates} In \cite{AshBen6}, Ashbaugh and Benguria proved, using
a wide variety of methods, the following upper bounds for $y^\ast(x)$:
\begin{equation*} \label{eq:ABest1}
  y^\ast(x)<K_2x\qquad\text{for}\quad  1<x\leq 1.396^{-},
\tag{AB${}_1$}
\end{equation*}
\begin{equation*} \label{eq:ABest2}
  y^\ast(x)\leq 1+x+\sqrt{2x-(1+x^2)/2}\qquad\text{for}\quad
                                      1.396^{-}\leq x \leq 1.634^{-},
\tag{AB${}_2$}
\end{equation*}
\begin{equation*} \label{eq:ABest3}
  y^\ast(x) \leq F(x)\qquad\text{for}\quad  1.634^{-} \leq x \leq 1.676^{-},
\tag{AB${}_3$}
\end{equation*}
\begin{equation*} \label{eq:ABest4}
  y^\ast(x) \leq H(x)-x\qquad\text{for}\quad  1.676^{-} \leq x \leq 2.198^{+},
\tag{AB${}_4$}
\end{equation*}
and
\begin{equation*} \label{eq:ABest5}
  y^\ast(x) \leq G(x)\qquad\text{for}\quad  2.198^{+} \leq x \leq 2.539^{-},
\tag{AB${}_5$}
\end{equation*}
where the functions $H(x), F(x)$ and $G(x)$ are defined by
$$
H(x) =  \begin{cases}
        &6\qquad\text{for}\quad x=1\\
        &\min\limits_{1\leq\eta,\xi<x}\left(2\eta+\dfrac{4\beta(\beta+\gamma)^2
        (x-1)(x-\beta\gamma/(\beta+\gamma-1))^2}{(2\beta-1)(2\gamma-1)(x-\eta)(x-\xi)
        (4x-2-\eta-\xi)}\right)\\ 
	&\qquad\text{for}\quad x>1
  \end{cases}
$$
with $\beta = \eta+\sqrt{\eta^2 - \eta}$ and $\gamma = \xi + \sqrt{\xi^2 - \xi}$,
$F(x)$ is the middle root of the cubic
$$
  2xy^3-2(5x^2+3x+1)y^2+(6x^3+39x^2+2x-1)y-(24x^3+11x^2-4x-1)=0
$$
and
$$
  G(x) = \inf_{\beta>1/2}\left(\dfrac{\beta^2}{2\beta-1}+
         \dfrac{x-\beta^2/(2\beta-1)}{C_2(\beta)(x-\beta^2/(2\beta-1))-1}\right)
$$
with the infimum taken over values of $\beta$ satisfying
$x > \beta^2/(2\beta-1)+1/C_2(\beta)$ and with
$$
  C_2(\beta) = \dfrac{2\beta-1}{\beta}\dfrac{\dsp\int_0^{j_{0,1}}t^3J_0^{2\beta}(t)dt}
                                                {\dsp\int_0^{j_{0,1}}tJ_0^{2\beta}(t)dt}
$$
We remind that $J_0(t)$ denotes the standard Bessel function of order zero and $j_{0,1}$ is its
first positive zero. For the derivation of the bound \eqref{eq:ABest4} and more
discussion of it, see \cite{AshBen5}. The other bounds given above are due to Hong-Cang Yang (see
\cite{HCY} for \eqref{eq:ABest2}) and Ashbaugh and Benguria (see \cite{AshBen2},
\cite{AshBen4} for \eqref{eq:ABest1}, \cite{AshBen5} for \eqref{eq:ABest5}, and
\cite{AshBen6} for \eqref{eq:ABest3}).

\begin{figure}[t]
\begin{center}
\includegraphics[width=0.9\textwidth,clip=yes]{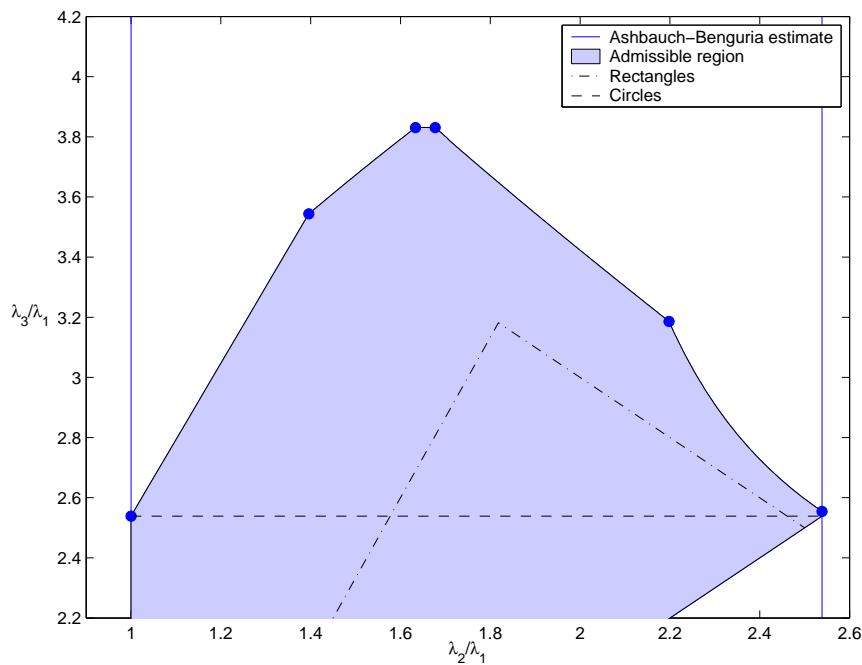}
\end{center}
\caption[]{\label{fig:figure_1}\small Admissible range (shaded) of $\dsp\left(\frac{\lam_2}{\lam_1},\frac{\lam_3}{\lam_1}\right)$ according to \cite{AshBen6}. Shown for comparison 
are the maximum values of $\dsp\frac{\lam_3}{\lam_1}$ as functions of
$\dsp\frac{\lam_2}{\lam_1}$ for rectangles and disjoint unions of circles.}
\end{figure}

The admissible region for $\dsp\left(\frac{\lam_2}{\lam_1},\frac{\lam_3}{\lam_1}\right)$ defined by
\eqref{eq:Ks}, obvious bounds $\dsp\frac{\lam_2}{\lam_1}\ge 1$ and
$\dsp\frac{\lam_3}{\lam_1}\ge \frac{\lam_2}{\lam_1}$, and inequalities \eqref{eq:ABest1}--\eqref{eq:ABest5},  is
shown in Figure~\ref{fig:figure_1}.

We make two remarks following  \cite{AshBen6}:
\begin{rem}
The inequalities \eqref{eq:ABest1}--\eqref{eq:ABest5} all apply on
broader intervals of $x$-values than the intervals specified explicitly with them; the given
intervals indicate the range for which the corresponding inequality gives the best bound
yet found.
\end{rem}

\begin{rem}
The absolute maximum of the right-hand sides of \eqref{eq:ABest1}--\eqref{eq:ABest5} occurs at the
point where $F(x)$ has a maximum within the interval where it is the best bound.  That happens
at the point $(x,y)\approx(1.65728,3.83103)$, and
imply the best upper bound \eqref{eq:lam3best} yet proven for $\dsp\frac{\lambda_3}{\lambda_1}$.
\end{rem}

\section{Numerical analysis of random domains}

To the best of our knowledge, there have been no large scale numerical experiments on low eigenvalues of the
Dirichlet Laplacian for planar domain. In an attempt to improve the existing estimates on the range of
$\dsp\left(\frac{\lam_2}{\lam_1},\frac{\lam_3}{\lam_1}\right)$, we conducted such experiments for a variety
of domain classes.

\subsection*{General method} For each particular domain, the calculation of the first three 
eigenvalues has
been conducted using a standard finite element method
implementation via PDEToolbox \cite{PDETool} and FEMLAB
\cite{FEMLAB} in Matlab, with two or three mesh refinements. For
simple domains with relatively ``high'' values of the ratio
$\dsp\frac{\lam_3}{\lam_1}$, optimization with respect to the
parameters describing domains of the particular class was
performed in order to maximize this ratio. The results of
calculations for some classes of domains are described below and
summarized at the end of this Section.

For each class of the domains we represent the results in the
following graphical form. The range $[1,K_2]$ of possible values of 
$x=\dsp\frac{\lam_2}{\lam_1}$ is split into subintervals of length $\del x$ 
(normally $\approx 0.05$). In each subinterval we choose, if it exists, a domain
with maximal $y=\dsp\frac{\lam_3}{\lam_1}$ and plot the corresponding point $(x,y)$.
For comparison, the graphs of $\dsp \left.y(x)=y^\ast(x)\right|_{\text{rectangles}}$ and/or
$\dsp \left.y^\ast(x)\right|_{\text{circles}}$ are shown.

\subsection*{Triangles, quadrilaterals and ellipses} 
\begin{figure}[thb!]
\begin{center}
\includegraphics[width=0.9\textwidth,clip=yes]{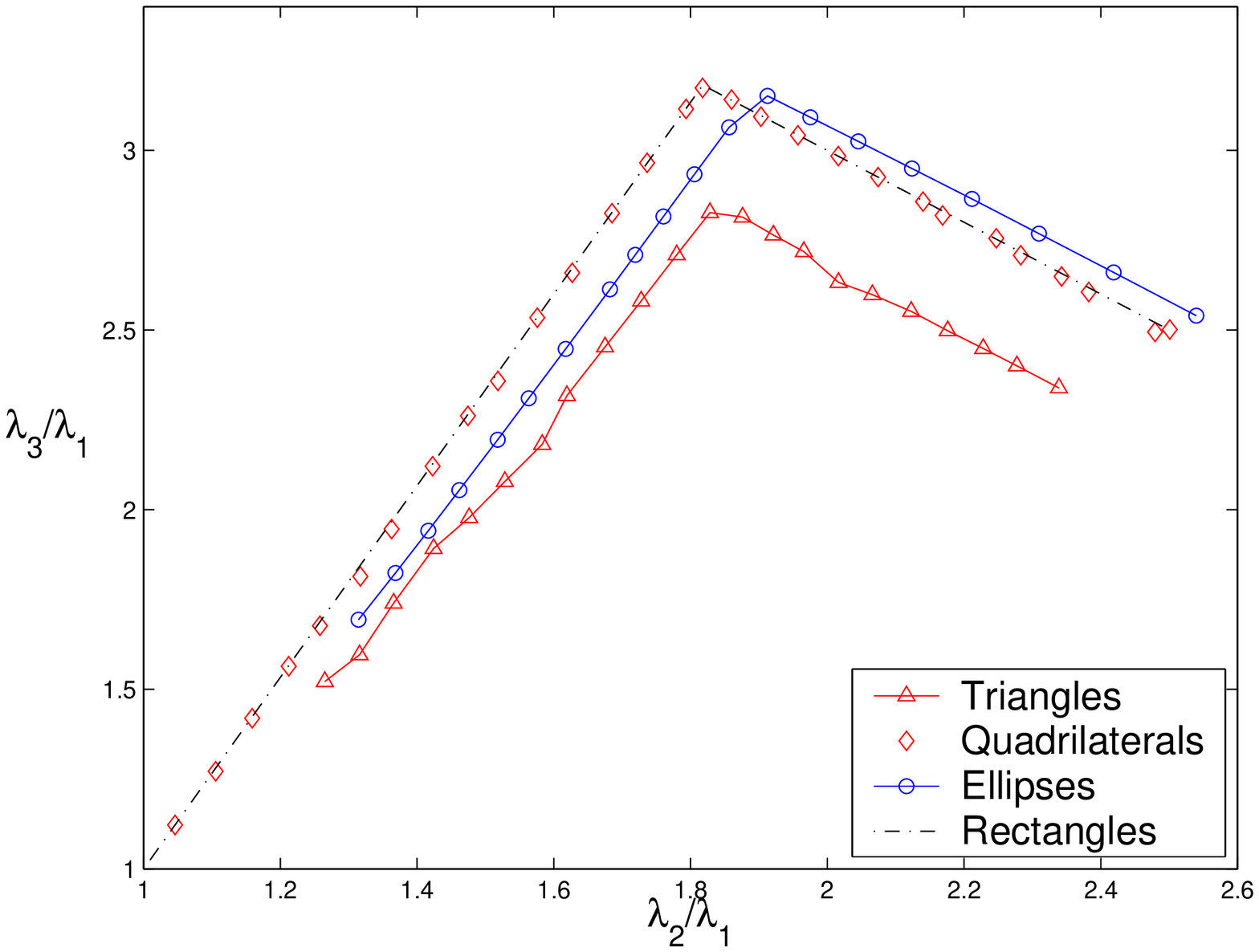}
\end{center}
\caption[]{\label{fig:tri_quad_ell}\small $y^\ast(x)$ for triangles, quadrilaterals and ellipses}
\end{figure}

We start with cyclic calculations for all triangles, 
with angle step $2.5^{\circ}$. For the triangles with relatively high ratio
of $\lambda_3/\lambda_1$ we repeat the procedure in the local
neighbourhood with angle step $0.5^{\circ}$. The results are
shown in Fig.~\ref{fig:tri_quad_ell}.

\begin{figure}[thb!]
\begin{center}
\includegraphics[width=0.9\textwidth,clip=yes]{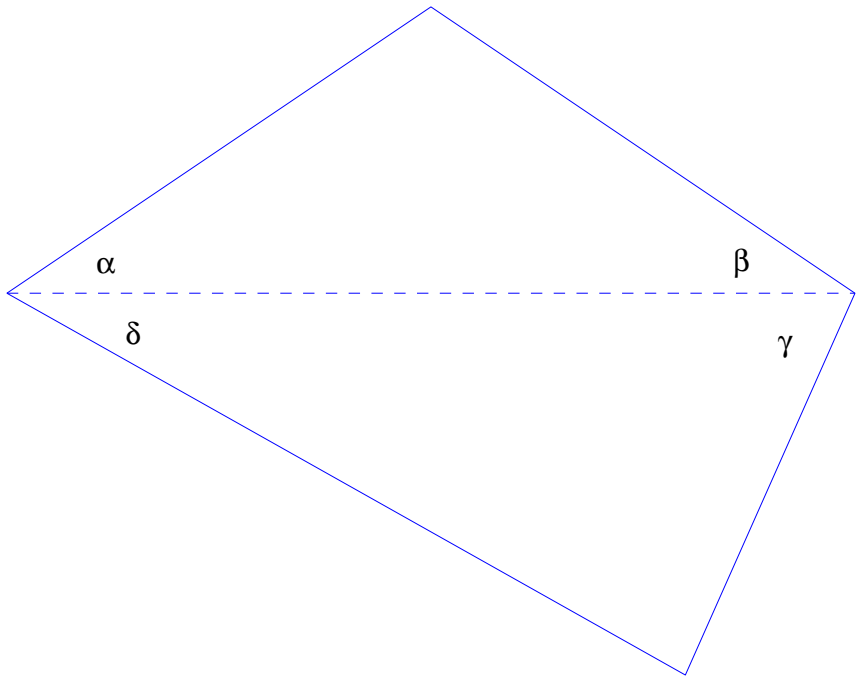}
\end{center}
\caption[]{\label{fig:quadrilat}Parametrization of a quadrilateral.}
\end{figure}

The computational procedure for quadrilaterals is essentially
the same as the one for triangles, with parameters $\alpha$,
$\beta$, $\gamma$ and $\delta$ in the region $(0,\pi)$ (see
Fig.~\ref{fig:quadrilat}).
We choose an angle step of $2.5^{\circ}$. Note
that quadrilaterals with negative $\alpha$ and $\beta$, or
$\gamma$ and $\delta$ do not have to be considered separately --- 
they fit into  the scheme above by choosing another diagonal as
a starting point and re-scaling. In case of relatively high ratio
$\lambda_3/\lambda_1$ ($\geq 3$) we repeated the calculation with
an angle step $0.5^{\circ}$ in the local neighbourhood of that
quadrilateral. The results are shown in
Fig.~\ref{fig:tri_quad_ell}.

The results of cyclic calculations for the ellipses, with axis ratio
varying between $1$ and $5$ with step $0.1$ ($0.001$ in the vicinity of the ellipse 
with highest $\lambda_3/\lambda_1$), are also shown in Fig.~\ref{fig:tri_quad_ell}.

\begin{rem}\label{rem:surprise}
Rather surprisingly, Fig.~\ref{fig:tri_quad_ell} suggests that 
\begin{equation}\label{eq:rect_quad}
\left.y^\ast(x)\right|_{\text{quadrilaterals}} \approx \left.y^\ast(x)\right|_{\text{rectangles}}
\end{equation}
We give a partial explanation of this fact in the next Section, see Remark~\ref{rem:rect_eq_quad}.
\end{rem}

\subsection*{Annuli and random sectors of annuli} The calculations for annuli with with inner radius $1$ and
outer radius $r$ demonstrated that the corresponding value $\dsp\frac{\lam_3}{\lam_1}(r)$ is monotonically increasing 
from $1$ to $K_2$ as $r$ changes from $1$ to $\infty$ (although convergence, for large $r$, is very slow --- just
logarithmic). These results  are not very informative and we do not include them in the graphs or the summary table 
below.

In calculations for sectors of the annuli of angle $\theta$, we choose $r$ randomly in the interval $(1,20)$ and 
$\theta$ randomly in the interval $(0.01\pi, 1.99\pi)$.  The results of calculations are shown in 
Fig.~\ref{fig:ngons_sectors}

\subsection*{Pseudo-random polygons} For polygons with more than four vertices, cyclic calculations
through all possible values of the geometric parameters with some reasonable step become impractical
due to the time constraints. Instead, we choose to perform calculations for randomly generated 
polygons. We employ the following simple procedure for generating a pseudo-random polygon with
$N$ vertices $\vbold_1,\dots,\vbold_N$ lying inside a square $[0,1]^2$.

\begin{description}
\item[Vertices $\vbold_1,\vbold_2,\vbold_3$] are chosen randomly using any pseudo-random generator.
\item[Vertices $\vbold_j$, $j=4,\dots,N-1$.] We choose a possible vertex at random. If the interval
$[\vbold_{j_1},\vbold_j]$ intersects any of previously constructed sides  $[\vbold_{k-1},\vbold_k]$, 
$k=1,\dots,j-2$, then we make another random choice.
\item[Vertex $\vbold_N$] is constructed in the same manner, but we additionally check that the
interval $[\vbold_N,\vbold_1]$ does not intersect any of the existing sides.
\end{description}
To avoid infinite loops, we abort the construction if the number of attempts at some stage exceeds
some sufficiently big number (say, 200). We also put in place a restriction forbidding very small angles
(which require special efforts in mesh generation).

\begin{figure}[thb!]
\begin{center}
\includegraphics[width=0.9\textwidth,clip=yes]{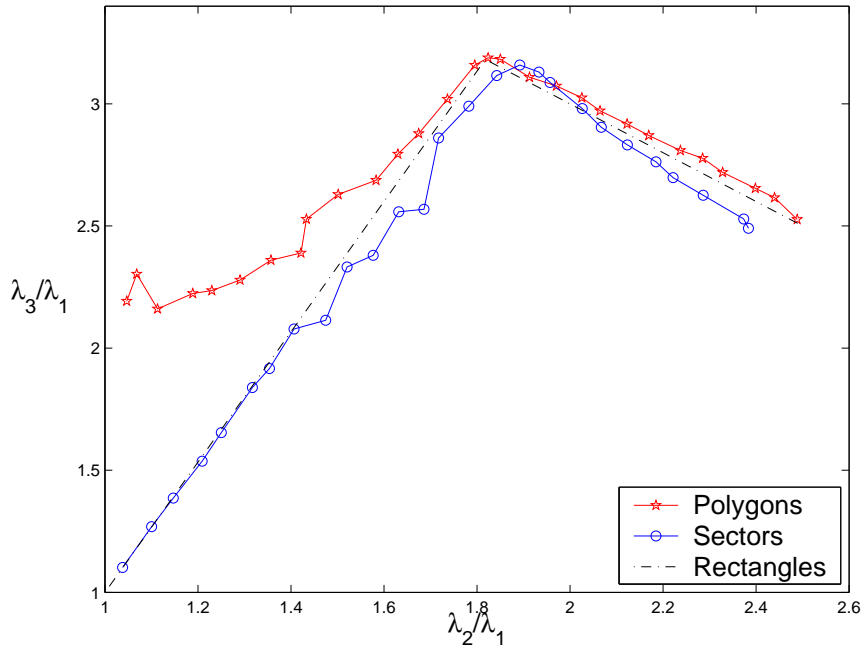}
\end{center}
\caption[]{\label{fig:ngons_sectors}\small $y^\ast(x)$ for random sectors of annuli and pseudo-random polygons}
\end{figure}

The collated results of calculations for pseudo-random pentagons, hex\-a\-gons and decagons are shown in Fig.~\ref{fig:ngons_sectors}. These results also include experiments on random perturbations of the 
rectangles constructed in the following way: $N$ points were randomly chosen on the sides of the rectangle $R_a$, with $a\in(1,5)$, and $1\le N\le 8$, and these points and the four vertices of the original rectangle
were randomly moved by a distance not exceeding $0.1a$ to form an $(N+4)$-gon.

\subsection*{Star-shaped domains (simply and non-simply connected)} The procedure described above for the polygons 
does not work very effectively for polygons with large number of vertices --- it often takes a long time to generate a
suitable vertex $\vbold_j$ with $j\gtrsim 10$.  Thus, in these cases we restrict ourselves to star-shaped
polygonal domains which are much easier to construct. Namely, for the vertices  $\dsp\vbold=re^{i\theta}$,
we choose the angles $\theta$ randomly between $0$ and $2\pi$, and the radii $r$ randomly between given
numbers $r_1$ and $r_2$. We conducted a series of experiments with a fixed number of vertices (13, 17 and
23), as well as a series of runs where the number of vertices was chosen randomly between four and thirty.

Additionally, we conducted a series of experiments of non-simply connected domains of the types
$\dsp R_{\sqrt\frac{8}{3}}\setminus S_1$, $\dsp S_2\setminus R_{\sqrt\frac{8}{3}}$, and $S_1\setminus S_2$,
where $S_j$ are random star-shaped polygons such that $\dsp S_1\subset  R_{\sqrt\frac{8}{3}}\subset S_2$,
and $R_{\sqrt\frac{8}{3}}$ is the rectangle with the maximum $\dsp\frac{\lam_3}{\lam_1}$.

\begin{figure}[thb!]
\begin{center}
\includegraphics[width=0.9\textwidth,clip=yes]{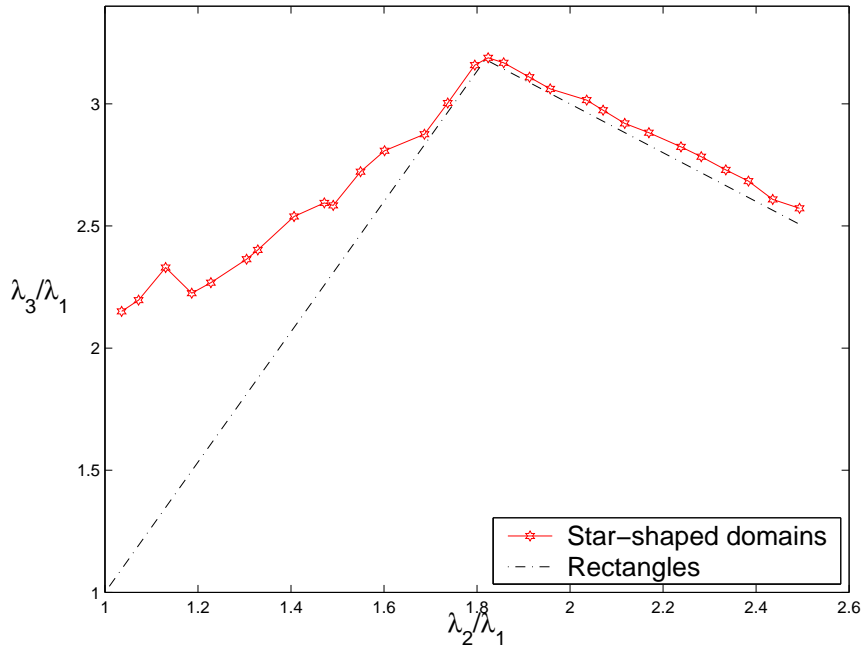}
\end{center}
\caption[]{\label{fig:star_shaped}\small $y^\ast(x)$ for pseudo-random star-shaped domains}
\end{figure}

The results for pseudo-random star-shaped domains are collated in Fig.~\ref{fig:star_shaped}.

\subsection*{Dumbbells and Jigsaw pieces} 
\begin{figure}[thb!]
\begin{center}
\includegraphics[width=0.9\textwidth,clip=yes]{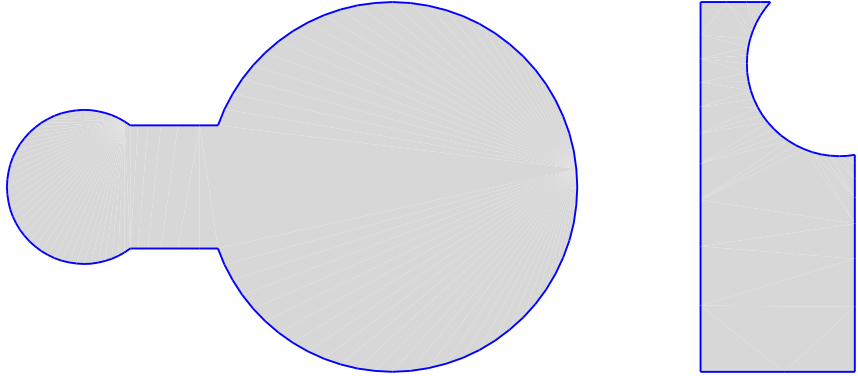}
\end{center}
\caption[]{\label{fig:dumb_jigsaw_typ}\small Typical dumbbell and jigsaw piece domains}
\end{figure}

By a {\it dumbbell\/} we understand a domain of the type 
\begin{equation}\label{eq:dumbbell}
([0,l]\times[-h,h])\cup C((0,0),r_1) \cup C((l,0),r_2)\,,
\end{equation}
where $l,h,r_1,r_2$ are positive parameters,
and $C(\vbold, r)$ denotes a circle with radius $r$ centred at $\vbold$. By a {\it jigsaw piece\/} we understand a domain 
of the type $R\setminus C$, where $R$ is a rectangle, and $C$ is a circle with a centre ``near'' the boundary
of the rectangle. Typical dumbbell and jigsaw piece domains are shown in Fig.~\ref{fig:dumb_jigsaw_typ}.

\begin{figure}[thb!]
\begin{center}
\includegraphics[width=0.9\textwidth,clip=yes]{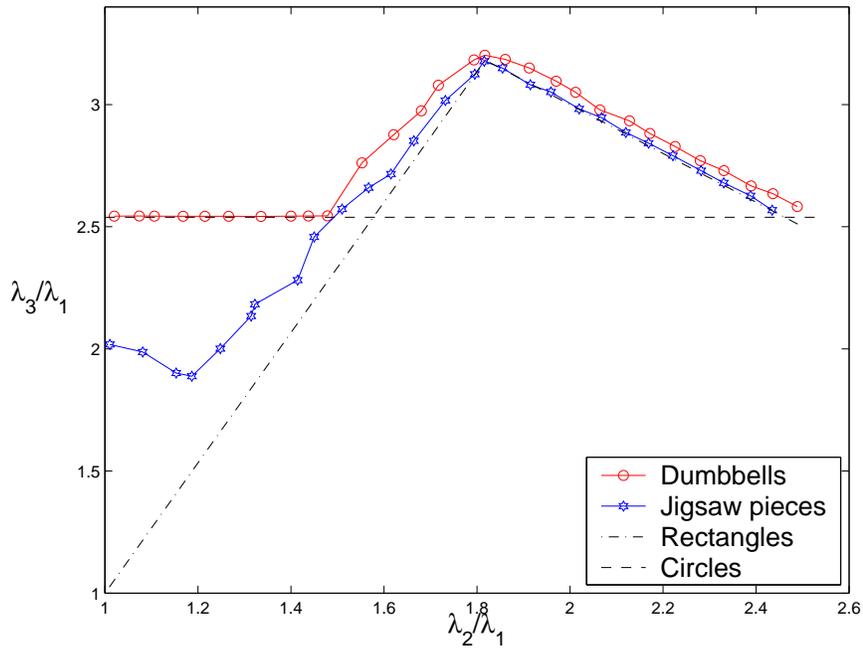}
\end{center}
\caption[]{\label{fig:dumbbells_jigsaws}\small $y^\ast(x)$ for pseudo-random dumbbells and jigsaw pieces}
\end{figure}

The results of numerical experiments on dumbbells and jigsaw pieces, with cyclical/random choice of the parameters, 
is shown in Fig.~\ref{fig:dumbbells_jigsaws}. For dumbbells, we also optimized over the parameters 
for domains with $\lam_3/\lam_1\approx 3.2$, allowing additionally the centres of the circles to move in the vertical direction along the sides of rectangles. However, this did not lead to the improvement of the results.

\subsection*{Summary of the numerical experiments} 

We summarize
the results of our numerical experiments in the following table.

\begin{table}[thb]
\begin{center}
\begin{tabular}{|l|r|r|r|}\hline
\bf{Type of domains}&\bf{No. of experiments}&$Y^\ast$&$\del_4$\\ \hline\hline 
Triangles&2145&2.827&0.016\\ \hline
Quadrilaterals&13222&3.183&$1.6\cdot 10^{-5}$\\ \hline
Sectors&360&3.149&$4.6\cdot 10^{-5}$\\ \hline
Ellipses&142&3.167&$9.1\cdot 10^{-5}$\\ \hline 
Random polygons&26867&3.189&$6.1\cdot 10^{-4}$\\ \hline 
Random star-shaped domains&18320&3.159&0.022\\ \hline 
Dumbbells&2871&3.202&$1.1\cdot 10^{-4}$\\ \hline
Jigsaw pieces&1420&3.178&0.003\\ \hline\hline
\bf{Total}&65337&3.202&$1.1\cdot 10^{-4}$\\ \hline
\end{tabular}
\end{center}
\caption[]{\label{tab:summarytable}\small Summary statistics for
numerical experiments. The value in the fourth column is
$\del_4:=\dsp \frac{\lam_4-\lam_3}{\lam_3}(\Ome^\ast)$, where
$\Ome^\ast$ is the domain which maximizes the ratio $\dsp
\frac{\lam_3}{\lam_1}=Y^\ast$ in the corresponding class of
domains}
\end{table}

As seen in the last column  of Table~\ref{tab:summarytable}, in each class of domains the maximum of the ratio
$\dsp\frac{\lam_3}{\lam_1}$ is attained, within the accuracy of computations, on a domain with {\it degenerate}
eigenvalue $\lam_3\approx \lam_4$. The same, of course, holds for rectangles, see \eqref{eq:yastrect}. This 
allows us to conjecture that the absolute maximum and any local maxima
of $\lam_3\approx \lam_4$ are also attained on domains with degenerate $\lam_3$. We give a partial proof of this
conjecture in the next Section.

The computed absolute maximum ratio $Y^\ast\approx 3.202$ is attained on the dumbbell-shaped domain \eqref{eq:dumbbell} 
with $l=1$, $h=1.4510$, $r_1=0.7814$, and $r_2=0.7818$, see Fig.~\ref{fig:optimal_dumbbell}. Note that the maximum value $Y^\ast$
is only slightly higher than the corresponding value $\dsp\left.Y^\ast\right|_{\text{rectangles}}\approx 3.1818$.

\begin{rem} Additional experiments were conducted in order to check whether a maximizer is likely to be a simply 
connected domain. Namely, for the dumbbell-shaped domain $\Ome$ described above, we computed the eigenvalues 
for a number of domains obtained by removing a small hole from $\Ome$. In all the cases, the ratio of the third and the
first eigenvalue for a perturbed problem was quite significantly less than that for $\Ome$.
\end{rem}

\begin{figure}[thb!]
\begin{center}
\includegraphics[width=0.9\textwidth]{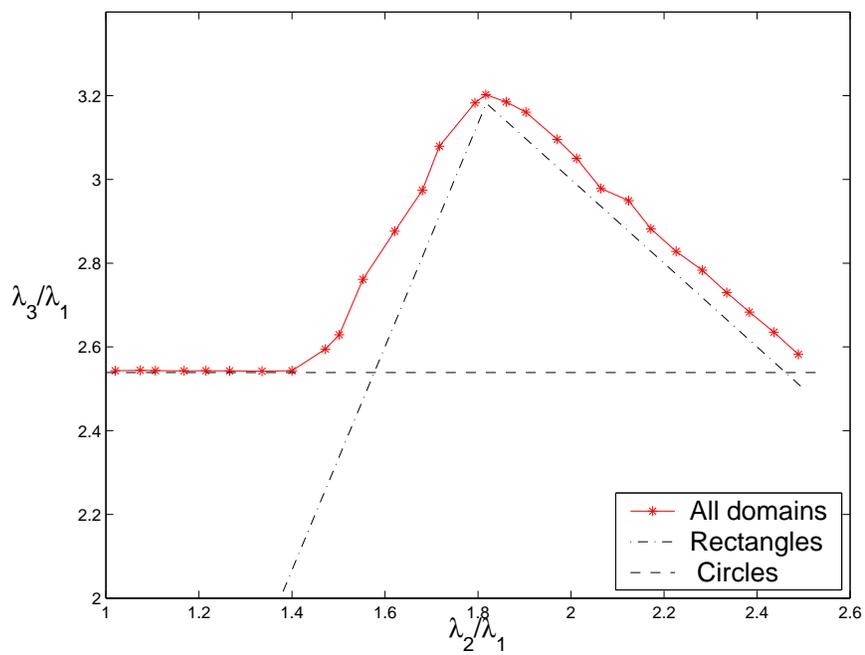}
\end{center}
\caption[]{\label{fig:finaldata}\small $y^\ast(x)$ for all computed domains}
\end{figure}

The graph of the function $y^\ast(x)$ built on the basis of all numerical experiments is shown in  
Fig.~\ref{fig:finaldata}.

\clearpage

\section{Asymptotic results}

\renewcommand{\thefootnote}{\fnsymbol{footnote}}

In this Section, using standard perturbation techniques, we establish several results which, although
don't give the full answer to the question of maximizing the ratio $\dsp\frac{\lam_3}{\lam_1}(\Ome)$ among all 
planar domains $\Ome$, give some indication which of the domains may or may not be a maximizer. 
We first proof the following
\begin{thm}\label{thm:rect} 
The rectangle $R_{\sqrt{\frac{8}{3}}}$ does not maximize the $\dsp\frac{\lam_3}{\lam_1}$
among all planar domains.\footnote{We were informed by Niculae Mandrache that he had independently obtained a similar result.}
\end{thm}

This should be compared, however, with Remark~\ref{rem:rect_eq_quad} below.

We also give a proof of the following more general result, which justifies the remark made at the end of last Section.
\begin{thm}\label{thm:degen}
Suppose that $\Omega\subseteq\Rbb[2]$ is a local maximizer of $\dsp\frac{\lam_3}{\lam_1}(\Ome)$ among planar domains with sufficiently smooth boundary. Then $\lam_3(\Ome)=\lam_4(\Ome)$.
\end{thm}

We should emphasize here that neither the statement not the proof (found below) of Theorem~\ref{thm:degen} is fully
rigorous. In the former, we do not discuss the requirements on the smoothness of the boundary and the concept of a local maximizer; we also do not prove that maximizers actually exist. In the latter, we rely on the following unproven, although very plausible, conjecture.

\begin{conj}\label{conj:nodal} Let $\lam_3$ be a simple eigenvalue of the Dirichlet Laplacian on a planar connected domain. 
Then not all nodal lines of the corresponding eigenfunction are closed.
\end{conj} 

Such a conjecture is not unreasonable since, in general, it is quite difficult to construct domains for which
even one nodal line of a low eigenfunction is closed, see \cite{HoHoNa}.

Before giving the proofs of Theorems~\ref{thm:rect} and \ref{thm:degen}, we recall, without proof, some classical results from the domain perturbation theory. The details can be found, 
e.g., in  \cite{Rel, SHSP}.

\subsection*{Domain perturbations} For simplicity, we restrict ourselves to domains in $\Rbb[2]$; all the results
stated here hold in any dimension.

Consider, for small values of real parameter $|\eps|$, a family of bounded domains $\Ome^\eps$ in $\Rbb[2]$ of variable $\tilde\xbold=(\tilde x_1,\tilde x_2)$,
which are transformed by the change of coordinates
\begin{equation}\label{eq:Sb}
\xbold=\tilde\xbold+\eps \Sbold(\xbold)
\end{equation}
into the domain $\Ome=\Ome^0$ in  $\Rbb[2]$ of variable $\xbold$. We assume that the boundary $\partial\Ome$ 
and the vector-function $\Sbold$ are sufficiently smooth.

Let $\nbold$ be the outer unit normal to $\partial\Ome$, and denote 
$$
f=\Sbold\cdot\nbold
$$
(in fact, $\eps f$ is a smooth function on $\partial\Ome$ which gives, up to the leading order for small $\eps$,
the normal distance between $\partial\Ome$ and $\partial\Ome^\eps$).

Denote by $\lam_1<\lam_2\le\dots\le\lam_j\le\dots$ the eigenvalues of the Dirichlet Laplacian on $\Ome$ and
by $\{u_j\}$ the corresponding basis of normalized orthogonal eigenfunctions (which are chosen real). Also, 
denote by $\lam_j^\eps$ the eigenvalues of the Dirichlet Laplacian on $\Ome^\eps$. For sufficiently small $|\eps|$, the  $\lam_j^\eps$ 
are continuous functions of $\eps$ and tend to $\lam_j$ as $\eps\to 0$.

The following two results go back to Rellich.

\begin{prop}\label{prop:per1} 
Let $\lam_j$, $j\ge 1$, be a {\em simple} eigenvalue of the  Dirichlet Laplacian on $\Ome$. Then $\lam_j^\eps$ has
the asymptotic expansion
\begin{equation}\label{eq:asympt}
\lam_j^\eps=\lam_j+\eps\tilde\lam_{j, 1}+\eps^2\tilde\lam_{j, 2}+\dots\,,
\end{equation}
as $\eps\to 0$, where
\begin{equation}\label{eq:cor1}
\tilde\lam_{j, 1}=-\int_{\partial\Ome} f\,\left|\frac{\partial u_j}{\partial n}\right|^2\,d\sig\,.
\end{equation}
\end{prop}

The situation is slightly more complicated when $\lam_j=\dots=\lam_{j+m}$ is an eigenvalue of multiplicity $m+1$.
For simplicity, we consider just the case $m=1$.

\begin{prop}\label{prop:per2} Let $\lam_k=\lam_{k+1}$ be a {\em double} eigenvalue of the Dirichlet Laplacian on $\Ome$. Then, as $\eps\to 0$, $\lam^\eps_k$ and $\lam^\eps_{k+1}$ still have the asymptotic expansions \eqref{eq:asympt} ($j=k, k+1$) with
\begin{equation}\label{eq:cor2}
\tilde\lam_{k, 1}=\frac{1}{\eps}\,\min(\eps\mu_1, \eps\mu_2)\,,\qquad
\tilde\lam_{k+1, 1}=\frac{1}{\eps}\,\max(\eps\mu_1, \eps\mu_2)\,,
\end{equation}
where $\mu_1$, $\mu_2$ are two real roots of the quadratic equation
\begin{equation}\label{eq:quad}
(F_{k,k}+\mu)(F_{k+1,k+1}+\mu)-F_{k,k+1}^2=0
\end{equation}
and
\begin{equation}\label{eq:Fs}
F_{p,q}=\int_{\partial\Ome} f\,\frac{\partial u_p}{\partial n}\,\frac{\partial u_q}{\partial n}\,d\sig\,.
\end{equation}
\end{prop}

We will be in fact interested in the asymptotic expansion of $\dsp\frac{\lam_j^\eps}{\lam_1^\eps}$, which follows 
from \eqref{eq:asympt}:
\begin{equation}\label{eq:ratio_asympt}
\frac{\lam_j^\eps}{\lam_1^\eps}=\frac{\lam_j}{\lam_1}+
\frac{\eps}{(\lam_1)^2}(\tilde\lam_{j, 1}\lam_1-\tilde\lam_{1, 1}\lam_j)+O(\eps^2)\,,
\end{equation}

\subsection*{Proof of Theorem \ref{thm:rect}} Let $\Ome=R_{\sqrt{\frac{8}{3}}}$ be the rectangle 
$\dsp\left\{(x_1,x_2):0<x_1<1,\ 0<x_2<\sqrt{\frac{8}{3}}\right\}$. We shall construct an explicit perturbation 
$\Ome^\eps$ using \eqref{eq:Sb} such that the first correction term in the asymptotic formula 
\eqref{eq:ratio_asympt} is positive for $\eps>0$, and therefore $\dsp\frac{\lam_3^\eps}{\lam_1^\eps}>\frac{\lam_3}{\lam_1}$ for sufficiently small positive $\eps$.

Let 
$$
\Ome^\eps=\left\{(x_1,x_2):0<x_1<1,\ 0<x_2<\sqrt{\frac{8}{3}}+\eps g(x_1)\right\}\,,
$$
where 
$$
g(x_1)=c_0+\sum_{l=0}^\infty \sqrt{2} c_l \cos(\pi l x_1)\,.
$$
We will choose the coefficients $c_l$ later. 

The corresponding function $f$ appearing in the asymptotic formulae above is
$$
f(x_1,x_2)=\begin{cases}
g(x_1)\,,\qquad&\text{if}\quad x_2=\sqrt{\frac{8}{3}}\,,\ 0\le x_1\le 1\,,\\
0\,,\qquad&\text{if}\quad (x_1, x_2)\in\partial\Ome\,,\ x_2\ne\sqrt{\frac{8}{3}}\,.
\end{cases}
$$
 
Note that we shall use \eqref{eq:cor1} for computing $\tilde\lam_{1,1}$ and \eqref{eq:cor2} for computing $\tilde\lam_{3,1}$ and $\tilde\lam_{4,1}$, since $\lam_3=\lam_4$ is a double eigenvalue of the unperturbed problem.
Elementary but tedious calculations show that the correction terms $\tilde\lam_{k,1}$, $k=1,2,3,4$, depend
only upon the parameters $c_j$ with $j=0,\dots,4$. For brevity, we omit the explicit expressions.

Let us choose the parameters $c_0,\dots,c_4$ in such a way that $\tilde\lam_{3,1}=\tilde\lam_{4,1}$ (i.e., 
$\lam_3^\eps$ remains a double eigenvalue up to the linear terms in $\eps$). This, by Proposition~\ref{prop:per2},
happens when $F_{3,3}=F_{4,4}$ and $F_{3,4}=0$, which in turn leads to the following conditions on coefficients
$c_j$:
\begin{equation}\label{eq:eqshift}
c_3=c_1\,,\qquad c_4=9c_2-8\sqrt{2}c_0\,.
\end{equation} 

Under conditions \eqref{eq:eqshift}, asymptotic formula \eqref{eq:ratio_asympt} simplifies dramatically, and becomes
$$
\frac{\lam_3^\eps}{\lam_1^\eps}-\frac{\lam_3}{\lam_1}=\frac{96\sqrt{3}}{121}(c_2-\sqrt{2}c_0)\eps+O(\eps^2)\,,
$$
and we can choose $c_0$ and $c_2$ in such a way that its right-hand side is positive for sufficiently small positive 
$\eps$. This proves Theorem~\ref{thm:rect}.

\begin{rem}\label{rem:rect_eq_quad} Let $\Ome=R_a$ be any rectangle and consider the perturbations $\Ome^\eps$ as above but
with function $f$ {\bf linear\/} in $x_1, x_2$ (and, naturally, $a$ replacing $\sqrt{\frac{8}{3}}$ throughout). Thus, we are 
considering  {\bf quadrilaterals\/} $\Ome^\eps$ which are small perturbations of the rectangle $R_a$. The same elementary calculations then imply that, for $\dsp x^\eps:=\frac{\lam_2}{\lam_1}(\Ome^\eps)$ and $\dsp y^\eps:=\frac{\lam_3}{\lam_1}(\Ome^\eps)$, we obtain, up to and inclusive of the terms of order $\eps$, that
$$
y^\eps=\left. y^\ast(x^\eps)\right|_{\text{rectangles}}\,,
$$
where $\left. y^\ast(x)\right|_{\text{rectangles}}$ is given by the right-hand side of \eqref{eq:yastrect}. In other words, 
up to the terms of order $\eps$ the rectangles are local maximizers among all quadrilaterals which are sufficently 
``close'' to them, cf. Remark \ref{rem:surprise}.
\end{rem}

\subsection*{Proof of Theorem \ref{thm:degen}} Suppose that $\Ome$ is a planar domain with sufficiently smooth 
boundary which locally maximizes the ratio $\frac{\lam_3}{\lam_1}$ in the following sense: for any sufficiently 
smooth perturbation $\Ome^\eps$ determined by \eqref{eq:Sb} we have
\begin{equation}\label{eq:ineq}
\frac{\lam_3^\eps}{\lam_1^\eps}\le \frac{\lam_3}{\lam_1}\,.
\end{equation}

Assume additionally that $\lam_3$ is a {\emph simple} eigenvalue of the  Dirichlet Laplacian in the unperturbed domain $\Ome$. We shall show that this assumption leads to the contradiction with Conjecture \ref{conj:nodal}.

Since both $\lam_1$ and $\lam_3$ are simple eigenvalues, the asymptotic formula \eqref{eq:ratio_asympt} becomes,
in accordance with Proposition~\ref{prop:per1},
$$
\frac{\lam_3^\eps}{\lam_1^\eps}-\frac{\lam_3}{\lam_1}=
\frac{\eps}{(\lam_1)^2}\int_{\partial\Ome} f\left(\lam_3\left|\frac{\partial u_1}{\partial n}\right|^2
-\lam_1\left|\frac{\partial u_3}{\partial n}\right|^2\right)\,d\sig+O(\eps^2)\,.
$$
Now, as $\eps$ can be chosen both positive and negative, \eqref{eq:ineq} can hold only if 
$$
\int_{\partial\Ome} f\left(\lam_3\left|\frac{\partial u_1}{\partial n}\right|^2
-\lam_1\left|\frac{\partial u_3}{\partial n}\right|^2\right)\,d\sig=0\,,
$$
and since $f$ is an arbitrary smooth function, this requires
$$
\lam_3\left|\frac{\partial u_1}{\partial n}\right|^2
=\lam_1\left|\frac{\partial u_3}{\partial n}\right|^2 
$$
everywhere on $\partial\Ome$. But the normal derivative of the first eigenfunction of the Dirichlet Laplacian is
non-zero everywhere on the boundary, so the last formula implies that the third eigenfunction has the same property, and therefore all its nodal lines are closed, in contradiction with Conjecture \ref{conj:nodal}. 

\section{Final Remarks}

\begin{figure}[thb!]
\begin{center}
\includegraphics[width=0.9\textwidth]{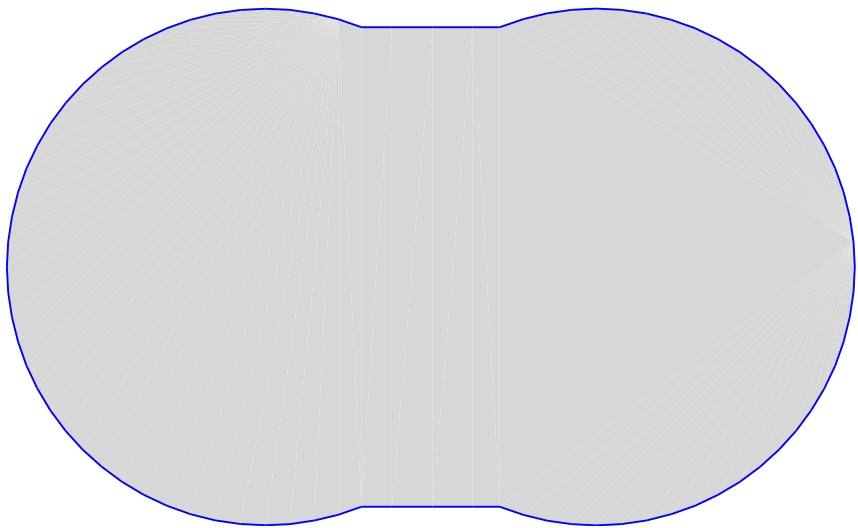}
\end{center}
\caption[]{\label{fig:optimal_dumbbell}\small Domain maximizing $\dsp \frac{\lam_3}{\lam_1}$ on the basis of computations}
\end{figure}

On the basis of the numerical computations and the results proven above we make the following
conjecture, most of which is still to be established (or disproved) rigorously:

{\it The domain maximizing the ratio $\dsp \frac{\lam_3}{\lam_1}$ for planar domains is close in shape to the optimal
computed dumbbell-shaped domain shown in Fig.~\ref{fig:optimal_dumbbell}, is simply-connected, and has a smooth boundary. The maximal admissible 
value $Y^\ast$ is approximately equal to or is slightly greater than $3.202$.}

\end{document}